\documentclass[10pt,oneside,reqno]{amsart}
\usepackage{eurosym}
\usepackage{latexsym}
\usepackage{amssymb,amsmath,amsfonts}
\usepackage[colorlinks=true]{hyperref}
\usepackage{graphicx}

\newtheorem{theorem}{Theorem}[section]
\newtheorem{lemma}[theorem]{Lemma}
\newtheorem{corollary}[theorem]{Corollary}
\newtheorem{proposition}[theorem]{Proposition}

\theoremstyle{definition}
\newtheorem{definition}[theorem]{Definition}

\newtheorem{example}[theorem]{Example}
\theoremstyle{remark}

\numberwithin{equation}{section}
\global\parskip 6pt \oddsidemargin
0.1cm \evensidemargin 0.4cm \headheight 0.1cm \topmargin -2.0cm
\begin{document}
\title{A new preconditioning algorithm for finding a zero of the sum of two
monotone operators and its application to image restoration problem}
\author{ $^{1}$ Ebru ALTIPARMAK AND $^{2}$ Ibrahim KARAHAN }
\address{$^{1}$ Department of Mathematics, Faculty of Science, Erzurum
Technical University, Erzurum, 25700, Turkey,}
\address{$^{2}$ Department of Mathematics, Faculty of Science, Erzurum
Technical University, Erzurum, 25700, Turkey,}
\email{$^{1}$ebru.altiparmak@erzurum.edu.tr,}
\email{$^{2}$ibrahimkarahan@erzurum.edu.tr,}

\begin{abstract}
Finding a zero of the sum of two monotone operators is one of the most
important problems in monotone operator theory, and the forward-backward
algorithm is the most prominent approach for solving this type of problem.
The aim of this paper is to present a new preconditioning forward-backward
algorithm to obtain the zero of the sum of two operators in which one is
maximal monoton and the other one is $M$-cocoercive, where $M$ is a linear
bounded operator. Furthermore, the strong convergence of the proposed
algorithm, which is a broader variant of previously known algorithms, has
been proven in Hilbert spaces. We also use our algorithm to tackle the
convex minimization problem and show that it outperforms existing
algorithms. Finally, we discuss several image restoration applications.
\end{abstract}

\maketitle

\section{Introduction}

Let $H$ be a real Hilbert space with inner product $\left\langle
.,.\right\rangle $ and the induced norm $\left\Vert .\right\Vert .$ One of
the most important problems in monotone operator theory is the problem of
finding a zero of the sum of two monotone operators so-called the monotone
inclusion problem which is defined by finding $x\in H$ such that
\begin{equation}
\text{ }0\in \left( A+B\right) \left( x\right)  \label{100}
\end{equation}%
where $A:H\rightarrow 2^{H}$ is a set-valued operator and $B:H\rightarrow H$
is an operator. This problem includes many mathematical problems such as
variational inequality problems, convex minimization problems, equilibrium
problems and convex-concave saddle point problems see e.g. : \cite%
{dixit,kitkuan,lopez,lorenz,sun,tseng}. More precisely, it has applications
in many scientific fields such as image processing, signal processing,
machine learning and statistical regression see e.g. : \cite{combet, sit,
sra, tib} . The most popular technique to solve the monotone inclusion
problem is the following forward-backward splitting algorithm which is
defined by Lions and Mercier \cite{lions}:
\begin{equation}
x_{n+1}=\left( I+\lambda _{n}A\right) ^{-1}\left( I-\lambda _{n}B\right)
x_{n},\text{ for all }n\in
\mathbb{N}
\label{101}
\end{equation}%
where $\lambda _{n}$ is a step size term and $A$ and $B$ are monotone
operators. If $B:H\rightarrow H$ is $1/L-$cocoercive operator and $\lambda
_{n}$ $\in \left( 0,2/L\right) ,$ the forward-backward splitting algorithm
converges weakly to a solution of the monotone inclusion problem. It is
well-known that the forward-backward splitting algorithm is a generalization
of classical proximal point and proximal gradient algorithm. Let $%
f:H\rightarrow
\mathbb{R}
$ be a differentiable convex function and let $g:H\rightarrow
\mathbb{R}
$ be a proper lower semi-continuous convex function$.$ The forward-backward
splitting algorithm (\ref{101}) is reduced to the proximal gradient
algorithm in this scenario, which is given as follows \cite{bau}:%
\begin{equation}
x_{n+1}=prox_{\lambda _{n}g}\left( I-\lambda _{n}\nabla f\right) x_{n},\text{
for all }n\in
\mathbb{N}
\label{102}
\end{equation}%
where $\lambda _{n}$ $>0$ is a step size. In subsequent work, Moudafi and
Oliny \cite{moudafi} introduced the following algorithm to solve to solve
the problem (\ref{100}) :%
\begin{equation}
\left\{
\begin{array}{l}
y_{n}=x_{n}+\theta _{n}\text{ }\left( x_{n}-x_{n-1}\right) \\
x_{n+1}=\left( I+\lambda _{n}A\right) ^{-1}\left( y_{n}-\lambda _{n}B\left(
x_{n}\right) \right) ,\text{ for all }n\in
\mathbb{N}%
\end{array}%
\right. ,  \label{106}
\end{equation}%
where $\theta _{n}$ is a inertial term on $\left[ 0,1\right)$. They studied
the weakly convergence of the algorithm, which satisfies the conditions $%
\sum_{n=1}^{\infty }\theta _{n}\left\Vert x_{n}-x_{n-1}\right\Vert
^{2}<\infty $ and $\lambda _{n}$ $<2/L$ where $L$ is the Lipschitz constant
of $B$. The presence of the inertial term increases the algorithm's
performance significantly.

In optimization problems, preconditioners are often used to speed up
first-order iterative optimization algorithms. For example, in gradient
descent method, one takes steps in the opposite direction of the gradient of
the function at the current point to find a local minimum of the real-valued
function. This algorithm is given by the following way:
\begin{equation*}
x_{n+1}=I-\lambda _{n}\nabla f\left( x_{n}\right) ,\text{for all }n\in
\mathbb{N}
.\text{ }
\end{equation*}%
The preconditioner $M $, which is a linear bounded operator, is applied to
the algorithm as follows:%
\begin{equation*}
x_{n+1}=I-\lambda _{n}M^{-1}\nabla f\left( x_{n}\right) ,\text{for all }n\in
\mathbb{N}
.
\end{equation*}%
The aim of the preconditioning is to change the geometry of the space to
make the level sets look like circles \cite{h}. In this situation, the
preconditioned gradient purposes getting closer to the extreme point and so
this accelerates the convergence. The classical splitting algorithms (\ref%
{101}) and (\ref{106}) may not generally be practical, and computing the
proximal mapping $\left( I+\lambda _{n}A\right) ^{-1}$ could be highly
costly. When we consider the preconditioned splitting algorithms with an
adequate mapping $M$, however, the algorithm becomes applicable.

In recent years, Lorenz and Pock \cite{lorenz} introduced the following
preconditioning algorithm to solve monotone inclusion problem:%
\begin{equation}
\left\{
\begin{array}{l}
y_{n}=x_{n}+\theta _{n}\text{ }\left( x_{n}-x_{n-1}\right) \\
x_{n+1}=\left( I+\lambda _{n}M^{-1}A\right) ^{-1}\left( I-\lambda
_{n}M^{-1}B\right) \left( y_{n}\right) ,\text{ for all }n\in
\mathbb{N}%
\end{array}%
\right. ,  \label{103}
\end{equation}%
where $\theta _{n}$ is an accelareted term on $\left[ 0,1\right) $ and $%
\lambda _{n}$ is a step size term. They proved the weak convergence of the
algorithm. It is clear that the Algorithm (\ref{103}) is reduced to the
classical forward-backward splitting algorithm (\ref{101}) for $\theta
_{n}=0 $ and $M=I$.

Subsequently, in 2021, Dixit et al. \cite{dixit} defined the following
algorithm which is called accelerated preconditioning forward-backward
normal $S-$iteration (APFBNSM):%
\begin{equation}
\left\{
\begin{array}{l}
y_{n}=x_{n}+\theta _{n}\text{ }\left( x_{n}-x_{n-1}\right) \\
x_{n+1}=\left( I+\lambda M^{-1}A\right) ^{-1}\left( I-\lambda M^{-1}B\right)
\left( \left( 1-\alpha _{n}\right) y_{n}+\alpha _{n}\left( I+\lambda
M^{-1}A\right) ^{-1}\left( I-\lambda M^{-1}B\right) \left( y_{n}\right)
\right) ,\text{ for all }n\in
\mathbb{N}%
\end{array}%
\right. ,  \label{104}
\end{equation}%
where $\alpha _{n}\in \left( 0,1\right) ,$ $\lambda \in \left[ 0,1\right) $
\ and $\theta _{n}\in \left[ 0,1\right) .$ They also proved weak convergence
of the proposed algorithm under some assumptions in a real Hilbert space $H$%
. For $\theta _{n}=0$ and $M=I,$ the accelerated preconditioning
forward-backward normal $S-$iteration (APFBNSM) is reduced to the normal $S-$%
iteration forward-backward splitting algorithm \cite{sahu} $\left(
nS-FBSA\right) :$%
\begin{equation*}
x_{n+1}=\left( I+\lambda A\right) ^{-1}\left( I-\lambda B\right) \left(
\left( 1-\alpha _{n}\right) y_{n}+\alpha _{n}\left( I+\lambda A\right)
^{-1}\left( I-\lambda B\right) \left( y_{n}\right) \right) ,\text{ for all }%
n\in
\mathbb{N}
.
\end{equation*}

In this paper, we present a new preconditioning forward-backward splitting
algorithm which generalizes many existed algorithms including the algorithms
(\ref{101}), (\ref{103}) and (\ref{104}), and which is more effective in
image restoration. Also, we prove that the sequence generated by the
proposed algorithm converges strongly to a solution of monotone inclusion
problem while the other algorithm's sequences converge weakly to the
solution of the same problem. The organization of this paper is listed as
follows. In the next section, we will give some definitions and lemmas to
study the convergence behaviour of the proposed algorithm. In Section 3, we
will present a new preconditioning forward-backward splitting algorithm and
study its convergence behaviour under mild restriction. In the last section,
we will give the application of the proposed algorithm to the image
restoration problem.

\section{ Preliminaries}

In this part, we will give some definitions and lemmas which play a
significant role in proving our main theorem. Let $C$ be a nonempty subset
of real Hilbert space $H$ and $T:C\rightarrow H$ be a mapping. A point $x\in
H$ is said to be a fixed point of $T$ if $Tx=x$ and the set of all fixed
point of $T$ is denoted by $F\left( T\right).$

\begin{definition}
\cite{cegi} Let $C$ be a nonempty subset of a real Hilbert space $H$ and $%
x\in H $. For any $z\in H$, if there exists a unique point $y\in C$ such
that
\begin{equation*}
\left\Vert y-x\right\Vert \leq \left\Vert z-x\right\Vert
\end{equation*}%
then $y$ is called the metric projection of $x$ onto $C$ and is denoted by $%
y=P_{C}x.$ If $P_{C}x$ exists and is uniquely determined for all $x\in H,$
then the operator $P_{C}:H\rightarrow C$ is called the metric projection.
\end{definition}

It is clear that the operator $P_{C}$ is nonexpansive and it can be
characterized by,%
\begin{equation*}
\left\langle x-P_{C}x,y-P_{C}x\right\rangle \leq 0\text{ for all }y\in C.
\end{equation*}

Let $A:H\rightarrow 2^{H}$ be a set-valued operator. If $\left\langle
u-v,x-y\right\rangle \geq 0$ for all $u\in Ax$ and $v\in Ay$, then $A$ is
said to be a monotone operator. If the graph of a monotone operator is not
properly contained in the graph of any other monotone operators, then $A$ is
said to be a maximal monotone operator.

Let $\ f:H\rightarrow \left( -\infty,+\infty \right] $ be a function. Then, $%
f$ is said to be proper if there exists at least one $x\in H$ such that $%
f\left( x\right)<+\infty $. Also, the subdifferential of a proper function $%
f $ is defined by%
\begin{equation*}
\partial f\left( x\right) =\left\{ u\in H:\left\langle y-x,u\right\rangle
\leq f\left( y\right) -f\left( x\right) \text{ for all }y\in H\right\} .
\end{equation*}%
and $f$ is subdifferentiable at $x\in H,$ if $\partial f\left( x\right) \neq
\emptyset .$ The elements of $\partial f\left( x\right) $ are called the
subgradients of $f$ at $x.$

\begin{definition}
\cite{bau} Let $\Gamma _{0}\left( H\right) $ denotes the class of all proper
lower semi-continuous convex functions defined from $H$ to $\left( -\infty
,+\infty \right].$ Let $g\in \Gamma _{0}\left( H\right) $ and $\phi >0.$ The
proximal operator of parameter $\phi $ of $g$ at $x$ is defined by
\begin{equation*}
prox_{\phi g}\left( x\right) =\arg \min_{y\in H}\left\{ g\left( y\right) +%
\frac{1}{2\phi }\left\Vert y-x\right\Vert ^{2}\right\} .
\end{equation*}
\end{definition}

\begin{example}
\cite{bau} \ Let $\phi \in \left( 0,+\infty \right) ,$ and let $x\in
\mathbb{R}
^n.$ Then, the proximal operator for $l_{1}-$norm is defined by
\begin{eqnarray*}
prox_{\phi \left\Vert .\right\Vert _{1}}\left( x\right) &=&\left( x-\phi
\right) _{+}-\left( -x-\phi \right) _{+} \\
&=&\left\{
\begin{array}{ccc}
x_{i}-\phi & if & x_{i}>\phi , \\
0 & if & -\phi \leq x_{i}\leq \phi , \\
x_{i}+\phi & if & x_{i}<-\phi ,%
\end{array}%
\text{ }\right.
\end{eqnarray*}
\end{example}

Let $M:H\rightarrow H$ be a bounded linear operator. $M$ is said to be
self-adjoint if $M^{\ast }=M$ where $M^{\ast }$ is the adjoint of operator $%
M $. A self-adjoint operator is said to be positive definite if $%
\left\langle M\left( x\right) ,x\right\rangle >0$ for every $0\neq x\in H$
\cite{limaye}. By using the self adjoint, positive and bounded linear
operator $M$, the $M$-inner product is defined by%
\begin{equation*}
\left\langle x,y\right\rangle _{M}=\left\langle x,M\left( y\right)
\right\rangle ,\text{ }\forall x,y\in H.
\end{equation*}%
In addition, the corresponding $M$-norm induced from the $M$-inner product
is defined by
\begin{equation*}
\left\Vert x\right\Vert _{M}^{2}=\left\langle x,M\left( x\right)
\right\rangle \text{ for all }x\in H.
\end{equation*}

\begin{definition}
\cite{dixit} Let $C$ be a nonempty subset of $H,$ $T:C\rightarrow H$ be an
operator and $M:H\rightarrow H$ be a positive definite operator. Then $T$ is
said to be:

\begin{enumerate}
\item[(i)] nonexpansive operator with respect to $M$-norm if%
\begin{equation*}
\left\Vert Tx-Ty\right\Vert _{M}\leq \left\Vert x-y\right\Vert _{M},\text{ }%
\forall x,y\in H,
\end{equation*}
\end{enumerate}
\end{definition}

\begin{enumerate}
\item[(ii)] $M$-cocoercive operator if$\left\Vert Tx-Ty\right\Vert
_{M^{-1}}^{2}\leq \left\langle x-y,Tx-Ty\right\rangle ,$ $\forall x,y\in H.$
\end{enumerate}

Similarly, $T$ said to be $k-$contraction mapping with respect to $M$-norm
if there exists $k\in \left[ 0,1\right) $ such that%
\begin{equation*}
\left\Vert Tx-Ty\right\Vert _{M}\leq k\left\Vert x-y\right\Vert _{M},\text{ }%
\forall x,y\in H.
\end{equation*}

\begin{proposition}
\cite{dixit} Let $A:H\rightarrow 2^{H}$ be a maximal monotone operator, $%
B:H\rightarrow H$ be a $M-$cocoercive operator, $M:H\rightarrow H$ be a
bounded linear self-adjoint and positive definite operator and $\lambda \in
\left( 0,1\right] $. Then we have the following properties:

\begin{enumerate}
\item[(i)] $I-\lambda M^{-1}B$ is nonexpansive with respect to $M$-norm,

\item[(ii)] $\left( I+\lambda M^{-1}A\right) ^{-1}$is nonexpansive with
respect to $M$-norm,

\item[(iii)] $J_{\lambda ,M}^{A,B}=\left( I+\lambda M^{-1}A\right)
^{-1}\left( I-\lambda M^{-1}B\right) $ is nonexpansive with respect to $M-$%
norm.
\end{enumerate}
\end{proposition}

\begin{proposition}
\cite{dixit} Let $A:H\rightarrow 2^{H}$ be a maximal monotone operator, $%
B:H\rightarrow H$ be a $M-$cocoercive operator, $M:H\rightarrow H$ be a
linear bounded self-adjoint and positive definite operator and $\lambda \in
\left( 0,\infty \right) $. Then $x\in H$ is a solution of monoton inclusion
problem (\ref{100}) if and only if $x$ is a fixed point of $J_{\lambda
,M}^{A,B}.$
\end{proposition}

\begin{lemma}
\cite{go}\label{dort} Let $C$ be a nonempty closed and convex subset of a
real Hilbert space $H$ and let $T:C\rightarrow H$ be a nonexpansive operator
with $F\left( T\right) \neq \emptyset .$ Then the mapping $I-T$ is
demiclosed at zero, that is, for any sequences $\left\{ x_{n}\right\} \in $ $%
H$ such that $x_{n}\rightharpoonup x$ $\in H$ and $\left\Vert
x_{n}-Tx_{n}\right\Vert \rightarrow 0$ as $n\rightarrow \infty ,$ then it
implies $x\in F\left( T\right) .$
\end{lemma}

\begin{lemma}
\cite{bau}\label{bir} Let $H$ be a real Hilbert space. Then for all \ $%
x,y\in H$ and $\lambda \in \left[ 0,1\right] ,$ the following properties are
hold:

\begin{enumerate}
\item[(i)] $\left\Vert x\pm y\right\Vert ^{2}=\left\Vert x\right\Vert
^{2}\pm 2\left\langle x,y\right\rangle +\left\Vert y\right\Vert ^{2},$

\item[(ii)] $\left\Vert x+y\right\Vert ^{2}\leq \left\Vert x\right\Vert
^{2}+2\left\langle y,x+y\right\rangle ,$

\item[(iii)] $\left\Vert \lambda x+\left( 1-\lambda \right) y\right\Vert
^{2}=\lambda \left\Vert x\right\Vert ^{2}+\left( 1-\lambda \right)
\left\Vert y\right\Vert ^{2}-\lambda \left( 1-\lambda \right) \left\Vert
x-y\right\Vert ^{2}.$
\end{enumerate}
\end{lemma}

\begin{lemma}
\cite{xu}\label{iki} Let $\left\{ s_{n}\right\} $ and $\left\{ \varepsilon
_{n}\right\} $ be sequences of nonnegative real numbers such that
\begin{equation*}
s_{n+1}\leq \left( 1-\delta _{n}\right) s_{n}+\delta _{n}t_{n}+\varepsilon
_{n},
\end{equation*}%
where $\left\{ \delta _{n}\right\} $ is a sequence in $\left[ 0,1\right] $
and $\left\{ t_{n}\right\} $ is a real sequence. If the following conditions
are hold, then $\lim_{n\rightarrow \infty }s_{n}=0:$

\begin{enumerate}
\item[(i)] $\sum\limits_{n=1}^{\infty }\delta _{n}=\infty ,$

\item[(ii)] $\sum\limits_{n=1}^{\infty }\varepsilon _{n}<\infty ,$

\item[(iii)] $\limsup_{n\rightarrow \infty }t_{n}\leq 0$
\end{enumerate}
\end{lemma}

\begin{lemma}
\cite{ma}\label{uc} Let $\left\{ \Phi _{n}\right\} $ be a sequence of real
numbers that does not decrease at infinity such that there exists a
subsequence $\left\{ \Phi _{n_{i}}\right\} $ of $\left\{ \Phi _{n}\right\} $
which satisfies $\Phi _{n_{i}}<\Phi _{n_{i+1}}$ for all $i\in
\mathbb{N}
.$ Let $\left\{ \tau \left( n\right) \right\} _{n\geq n_{0}}$ be a sequence
of integer which defined by:%
\begin{equation*}
\tau \left( n\right) :=\max \left\{ l\leq n:\Phi _{l}<\Phi _{l+1}\right\} .
\end{equation*}%
Then the following are satisfied:

\begin{enumerate}
\item[(i)] $\tau \left( n_{0}\right) \leq \tau \left( n_{0}+1\right) $ $\leq
...$ and $\tau \left( n\right) \rightarrow \infty ,$

\item[(ii)] $\Phi _{\tau \left( n\right) }\leq \Phi _{\tau \left( n\right) +1%
\text{ }}$and $\Phi _{n}\leq \Phi _{\tau \left( n\right) +1\text{ }},$ for
all $n\geq n_{0}.$
\end{enumerate}
\end{lemma}

\noindent

\section{Main Results}

In this section, we define a new preconditioning forward-backward splitting
algorithm and prove its strong convergence in real Hilbert space.

\begin{theorem}
\label{theorem} Let $M:H\rightarrow H$ be a bounded linear self-adjoint and
positive definite operator, $A:H\rightarrow 2^{H}$ be a maximal monotone
operator and $B:H\rightarrow H$ be a $M$-cocoercive operator such that $%
\Omega =\left( A+B\right) ^{-1}\left( 0\right) $ is nonempty. Let $f$ be a $%
k $-contraction mapping on $H$ with respect to $M$-norm and let $\lambda \in
\left( 0,1\right] .$ Let $\left\{ x_{n}\right\} $ be a sequence generated by%
\begin{equation}
\left\{
\begin{array}{l}
x_{0},x_{1}\in H \\
y_{n}=x_{n}+\theta _{n}\text{ }\left( x_{n}-x_{n-1}\right) \\
z_{n}=\left( I+\lambda M^{-1}A\right) ^{-1}\left( I-\lambda M^{-1}B\right)
\left( \left( 1-\alpha _{n}\right) y_{n}+\alpha _{n}\left( I+\lambda
M^{-1}A\right) ^{-1}\left( I-\lambda M^{-1}B\right) \left( y_{n}\right)
\right) \\
x_{n+1}=\beta _{n}f\left( z_{n}\right) +\left( 1-\beta _{n}\right) \left(
I+\lambda M^{-1}A\right) ^{-1}\left( I-\lambda M^{-1}B\right) \left(
z_{n}\right)%
\end{array}%
\right.  \label{105}
\end{equation}%
where $\left\{ \theta _{n}\right\} \subset \left[ 0,\theta \right] $ is a
sequence with $\theta \in \left[ 0,1\right) $ and $\left\{ \alpha
_{n}\right\} ,\left\{ \beta _{n}\right\} \in \left( 0,1\right) $ such that
the following conditions are hold:

\begin{enumerate}
\item[(i)] $0<a\leq \alpha _{n}\leq b<1$ for some $a,b\in
\mathbb{R}
,$

\item[(ii)] $0<c\leq \beta _{n}\leq d<1$ for some $c,d\in
\mathbb{R}
,$

\item[(iii)] $\sum_{n=1}^{\infty }\theta _{n}\left\Vert
x_{n}-x_{n-1}\right\Vert _{M}<\infty ,$

\item[(iv)] $\lim_{n\rightarrow \infty }\beta _{n}=0,$ $\sum_{n=1}^{\infty
}\beta _{n}=\infty .$
\end{enumerate}

Then the sequence $\left\{ x_{n}\right\} $converges strongly to a point $p$
in $\Omega $ where $p=P_{\Omega }f\left( p\right) .$

\begin{proof}
We will obtain the proof by dividing it into the following steps.

\textbf{Step 1 : } In this step, we show that the sequence $\left\{
x_{n}\right\} $ is bounded. Let $p\in \Omega $ such that $p=P_{\Omega
}f\left( p\right) $. Since $J_{\lambda ,M}^{A,B}$ is nonexpansive with
respect to $M$-norm, we obtain the followings from algorithm (\ref{105}):%
\begin{eqnarray}
\left\Vert y_{n}-p\right\Vert _{M} &=&\left\Vert x_{n}+\theta _{n}\text{ }%
\left( x_{n}-x_{n-1}\right) -p\right\Vert _{M}  \notag \\
&\leq &\left\Vert x_{n}-p\right\Vert _{M}+\theta _{n}\left\Vert \text{ }%
x_{n}-x_{n-1}\right\Vert _{M}  \label{1}
\end{eqnarray}%
and%
\begin{eqnarray}
\left\Vert z_{n}-p\right\Vert _{M} &=&\left\Vert J_{\lambda ,M}^{A,B}\left(
\left( 1-\alpha _{n}\right) y_{n}+\alpha _{n}J_{\lambda ,M}^{A,B}\left(
y_{n}\right) \right) -p\right\Vert _{M}  \notag \\
&\leq &\left\Vert \left( 1-\alpha _{n}\right) y_{n}+\alpha _{n}J_{\lambda
,M}^{A,B}\left( y_{n}\right) -p\right\Vert _{M}  \notag \\
&=&\left\Vert \left( 1-\alpha _{n}\right) \left( y_{n}-p\right) +\alpha
_{n}\left( J_{\lambda ,M}^{A,B}\left( y_{n}\right) -p\right) \right\Vert _{M}
\notag \\
&\leq &\left( 1-\alpha _{n}\right) \left\Vert y_{n}-p\right\Vert _{M}+\alpha
_{n}\left\Vert J_{\lambda ,M}^{A,B}\left( y_{n}\right) -p\right\Vert _{M}
\notag \\
&\leq &\left\Vert y_{n}-p\right\Vert _{M}.  \label{2}
\end{eqnarray}

Since $f$ is $k$-contractive mapping with respect to $M$-norm, we also
obtain the followings by combining (\ref{1}) and (\ref{2}):
\begin{eqnarray}
\left\Vert x_{n+1}-p\right\Vert _{M} &=&\left\Vert \beta _{n}f\left(
z_{n}\right) +\left( 1-\beta _{n}\right) J_{\lambda ,M}^{A,B}\left(
z_{n}\right) -p\right\Vert _{M}  \notag \\
&\leq &\left\Vert \beta _{n}\left( f\left( z_{n}\right) -p-f\left( p\right)
+f\left( p\right) \right) +\left( 1-\beta _{n}\right) \left( J_{\lambda
,M}^{A,B}\left( z_{n}\right) -p\right) \right\Vert _{M}  \notag \\
&\leq &\beta _{n}\left\Vert f\left( z_{n}\right) -f\left( p\right)
\right\Vert _{M}+\beta _{n}\left\Vert f\left( p\right) -p\right\Vert
_{M}+\left( 1-\beta _{n}\right) \left\Vert J_{\lambda ,M}^{A,B}\left(
z_{n}\right) -p\right\Vert _{M}  \notag \\
&\leq &\beta _{n}k\left\Vert z_{n}-p\right\Vert _{M}+\beta _{n}\left\Vert
f\left( p\right) -p\right\Vert _{M}+\left( 1-\beta _{n}\right) \left\Vert
z_{n}-p\right\Vert _{M}  \notag \\
&=&\left( 1-\beta _{n}\left( 1-k\right) \right) \left\Vert
z_{n}-p\right\Vert _{M}+\beta _{n}\left\Vert f\left( p\right) -p\right\Vert
_{M}  \notag \\
&\leq &\left( 1-\beta _{n}\left( 1-k\right) \right) \left\Vert
x_{n}-p\right\Vert _{M}+\beta _{n}.\frac{\theta _{n}}{\beta _{n}}\left\Vert
x_{n}-x_{n-1}\text{ }\right\Vert _{M}+\beta _{n}\left\Vert f\left( p\right)
-p\right\Vert _{M}.  \label{3}
\end{eqnarray}%
From the conditions $\left( ii\right) $ and $\left( iii\right) $, we have $%
\lim_{n\rightarrow \infty }$ $\frac{\theta _{n}}{\beta _{n}}\left\Vert
x_{n}-x_{n-1}\text{ }\right\Vert _{M}=0.$ So, there exists a positive
constant $K_{1}>0$ such that $\frac{\theta _{n}}{\beta _{n}}\left\Vert
x_{n}-x_{n-1}\text{ }\right\Vert _{M}\leq K_{1}.$ It follows from (\ref{3})
that,%
\begin{eqnarray*}
\left\Vert x_{n+1}-p\right\Vert _{M} &\leq &\left( 1-\beta _{n}\left(
1-k\right) \right) \left\Vert x_{n}-p\right\Vert _{M}+\beta _{n}\left(
K_{1}+\left\Vert f\left( p\right) -p\right\Vert _{M}\right) \\
&=&\left( 1-\beta _{n}\left( 1-k\right) \right) \left\Vert
x_{n}-p\right\Vert _{M}+\beta _{n}\left( 1-k\right) \left( \frac{%
K_{1}+\left\Vert f\left( p\right) -p\right\Vert _{M}}{\left( 1-k\right) }%
\right) \\
&\leq &\max \left\{ \left\Vert x_{n}-p\right\Vert _{M},\frac{%
K_{1}+\left\Vert f\left( p\right) -p\right\Vert _{M}}{\left( 1-k\right) }%
\right\} \\
&&\vdots \\
&\leq &\max \left\{ \left\Vert x_{1}-p\right\Vert _{M},\frac{%
K_{1}+\left\Vert f\left( p\right) -p\right\Vert _{M}}{\left( 1-k\right) }%
\right\}
\end{eqnarray*}%
for all $n\geq 1.$ This means that $\left\{ x_{n}\right\} $ is bounded so $%
\left\{ y_{n}\right\} ,\left\{ z_{n}\right\} $ are also bounded$.$

\textbf{Step 2 : }Next,\textbf{\ }we have to show that $x_{n}\rightarrow
p=P_{\Omega }f\left( p\right) .$ Indeed, using Lemma \ref{bir} we find the
followings for all $n\geq 1$:
\begin{eqnarray}
\left\Vert y_{n}-p\right\Vert _{M}^{2} &=&\left\Vert x_{n}+\theta _{n}\text{
}\left( x_{n}-x_{n-1}\right) -p\right\Vert _{M}^{2}  \notag \\
&\leq &\left\Vert x_{n}-p\right\Vert _{M}^{2}+2\theta _{n}\left\Vert
x_{n}-p\right\Vert _{M}\left\Vert x_{n}-x_{n-1}\right\Vert _{M}+\theta
_{n}^{2}\left\Vert x_{n}-x_{n-1}\right\Vert _{M}^{2}  \label{5}
\end{eqnarray}
and%
\begin{eqnarray}
\left\Vert z_{n}-p\right\Vert _{M}^{2} &=&\left\Vert J_{\lambda
,M}^{A,B}\left( \left( 1-\alpha _{n}\right) y_{n}+\alpha _{n}J_{\lambda
,M}^{A,B}\left( y_{n}\right) \right) -p\right\Vert _{M}^{2}  \notag \\
&\leq &\left\Vert \left( 1-\alpha _{n}\right) y_{n}+\alpha _{n}J_{\lambda
,M}^{A,B}\left( y_{n}\right) -p\right\Vert _{M}^{2}  \notag \\
&=&\alpha _{n}\left\Vert J_{\lambda ,M}^{A,B}\left( y_{n}\right)
-p\right\Vert _{M}^{2}+\left( 1-\alpha _{n}\right) \left\Vert
y_{n}-p\right\Vert _{M}^{2}-\alpha _{n}\left( 1-\alpha _{n}\right)
\left\Vert J_{\lambda ,M}^{A,B}\left( y_{n}\right) -y_{n}\right\Vert _{M}^{2}
\notag \\
&\leq &\left\Vert y_{n}-p\right\Vert _{M}^{2}-\alpha _{n}\left( 1-\alpha
_{n}\right) \left\Vert J_{\lambda ,M}^{A,B}\left( y_{n}\right)
-y_{n}\right\Vert _{M}^{2}  \notag \\
&\leq &\left\Vert y_{n}-p\right\Vert _{M}^{2}.  \label{4}
\end{eqnarray}
It follows from (\ref{5}), (\ref{4}), and Lemma \ref{bir} that%
\begin{eqnarray}
\left\Vert x_{n+1}-p\right\Vert _{M}^{2} &=&\left\Vert \beta _{n}f\left(
z_{n}\right) +\left( 1-\beta _{n}\right) J_{\lambda ,M}^{A,B}\left(
z_{n}\right) -p\right\Vert _{M}^{2}  \notag \\
&\leq &\left\Vert \beta _{n}\left( f\left( z_{n}\right) -f\left( p\right)
\right) +\left( 1-\beta _{n}\right) \left( J_{\lambda ,M}^{A,B}\left(
z_{n}\right) -p\right) +\beta _{n}\left( f\left( p\right) -p\right)
\right\Vert _{M}^{2}  \notag \\
&\leq &\left\Vert \beta _{n}\left( f\left( z_{n}\right) -f\left( p\right)
\right) +\left( 1-\beta _{n}\right) \left( J_{\lambda ,M}^{A,B}\left(
z_{n}\right) -p\right) \right\Vert _{M}^{2}+2\beta _{n}\left\langle f\left(
p\right) -p,x_{n+1}-p\right\rangle _{M}  \notag \\
&\leq &\beta _{n}\left\Vert f\left( z_{n}\right) -f\left( p\right)
\right\Vert _{M}^{2}+\left( 1-\beta _{n}\right) \left\Vert J_{\lambda
,M}^{A,B}\left( z_{n}\right) -p\right\Vert _{M}^{2}+2\beta _{n}\left\langle
f\left( p\right) -p,x_{n+1}-p\right\rangle _{M}  \notag \\
&\leq &\beta _{n}k^{2}\left\Vert z_{n}-p\right\Vert _{M}^{2}+\left( 1-\beta
_{n}\right) \left\Vert z_{n}-p\right\Vert _{M}^{2}+2\beta _{n}\left\langle
f\left( p\right) -p,x_{n+1}-p\right\rangle _{M}  \notag \\
&\leq &\left( 1-\beta _{n}\left( 1-k^{2}\right) \right) \left\Vert
z_{n}-p\right\Vert _{M}^{2}+2\beta _{n}\left\langle f\left( p\right)
-p,x_{n+1}-p\right\rangle _{M}  \notag \\
&\leq &\left( 1-\beta _{n}\left( 1-k^{2}\right) \right) \left[ \left\Vert
x_{n}-p\right\Vert _{M}^{2}+2\theta _{n}\left\Vert x_{n}-p\right\Vert
_{M}\left\Vert \text{ }x_{n}-x_{n-1}\right\Vert _{M}\right.  \notag \\
&&\left. +\theta _{n}^{2}\left\Vert \text{ }x_{n}-x_{n-1}\right\Vert _{M}^{2}%
\right] +2\beta _{n}\left\langle f\left( p\right) -p,x_{n+1}-p\right\rangle
_{M}  \notag \\
&\leq &\left( 1-\beta _{n}\left( 1-k^{2}\right) \right) \left\Vert
x_{n}-p\right\Vert _{M}^{2}+\theta _{n}\left\Vert \text{ }%
x_{n}-x_{n-1}\right\Vert _{M}\left[ 2\left\Vert x_{n}-p\right\Vert
_{M}\right.  \label{6} \\
&&\left. +\theta _{n}\left\Vert \text{ }x_{n}-x_{n-1}\right\Vert _{M}\right]
+2\beta _{n}\left\langle f\left( p\right) -p,x_{n+1}-p\right\rangle _{M}
\notag
\end{eqnarray}%
for all $n\geq 1.$ Since $\lim_{n\rightarrow \infty }$ $\theta
_{n}\left\Vert x_{n}-x_{n-1}\text{ }\right\Vert _{M}=0$, there exists a
positive constant $K_{2}>0$ such that $\theta _{n}\left\Vert x_{n}-x_{n-1}%
\text{ }\right\Vert _{M}\leq K_{2}$. From the inequality (\ref{6}) we
observe that, for all $n\geq 1,$
\begin{eqnarray}
\left\Vert x_{n+1}-p\right\Vert _{M}^{2} &\leq &\left( 1-\beta _{n}\left(
1-k^{2}\right) \right) \left\Vert x_{n}-p\right\Vert _{M}^{2}+3K_{3}\theta
_{n}\left\Vert \text{ }x_{n}-x_{n-1}\right\Vert _{M}  \notag \\
&&+2\beta _{n}\left\langle f\left( p\right) -p,x_{n+1}-p\right\rangle _{M}
\notag \\
&=&\left( 1-\beta _{n}\left( 1-k^{2}\right) \right) \left\Vert
x_{n}-p\right\Vert _{M}^{2}+3K_{3}\theta _{n}\left\Vert \text{ }%
x_{n}-x_{n-1}\right\Vert _{M}  \label{7} \\
&&+\beta _{n}\left( 1-k^{2}\right) \frac{2}{\left( 1-k^{2}\right) }%
\left\langle f\left( p\right) -p,x_{n+1}-p\right\rangle _{M},  \notag
\end{eqnarray}%
where $K_{3}=\sup_{n\geq 1}\left\{ \left\Vert x_{n}-p\right\Vert
_{M},K_{2}\right\} .$ In above inequality, if we take $\delta _{n}=\beta
_{n}\left( 1-k^{2}\right) ,$ $s_{n}=\left\Vert x_{n}-p\right\Vert _{M}^{2},$
$t_{n}=\frac{2}{\left( 1-k^{2}\right) }\left\langle f\left( p\right)
-p,x_{n+1}-p\right\rangle _{M}$ and $\varepsilon _{n}=3K_{3}\theta
_{n}\left\Vert \text{ }x_{n}-x_{n-1}\right\Vert _{M}$ then we have $%
s_{n+1}\leq \left( 1-\delta _{n}\right) s_{n}+\delta _{n}t_{n}+\varepsilon
_{n}$ for all $n\geq 1.$

Now, we want to show that $\limsup_{n\rightarrow \infty }\left\langle
f\left( p\right) -p,x_{n+1}-p\right\rangle _{M}\leq 0.$ So, we take into
account two cases to complete the proof.

First, we suppose that there exists $n_{0}\in
\mathbb{N}
$ such that $\left\{ \left\Vert x_{n}-p\right\Vert _{M}\right\} _{n\geq
n_{0}}$ is a nonincreasing sequence. So, the sequence $\left\{ \left\Vert
x_{n}-p\right\Vert _{M}\right\} $ is convergent since it is bounded from
below by $0.$ By using the condition $\left( iv\right) $, we have $%
\sum_{n=1}^{\infty }$ $\delta _{n}=\infty .$ We claim that $%
\limsup_{n\rightarrow \infty }\left\langle f\left( p\right)
-p,x_{n+1}-p\right\rangle _{M}\leq 0.$ By combining (\ref{5}) and (\ref{4})
with Lemma \ref{bir}, we get
\begin{eqnarray*}
\left\Vert x_{n+1}-p\right\Vert _{M}^{2} &=&\left\Vert \beta _{n}f\left(
z_{n}\right) +\left( 1-\beta _{n}\right) J_{\lambda ,M}^{A,B}\left(
z_{n}\right) -p\right\Vert _{M}^{2} \\
&=&\beta _{n}\left\Vert f\left( z_{n}\right) -p\right\Vert _{M}^{2}+\left(
1-\beta _{n}\right) \left\Vert J_{\lambda ,M}^{A,B}\left( z_{n}\right)
-p\right\Vert _{M}^{2}-\beta _{n}\left( 1-\beta _{n}\right) \left\Vert
f\left( z_{n}\right) -J_{\lambda ,M}^{A,B}\left( z_{n}\right) \right\Vert
_{M}^{2} \\
&\leq &\beta _{n}\left\Vert f\left( z_{n}\right) -p\right\Vert
_{M}^{2}+\left( 1-\beta _{n}\right) \left\Vert z_{n}-p\right\Vert _{M}^{2} \\
&\leq &\beta _{n}\left\Vert f\left( z_{n}\right) -p\right\Vert
_{M}^{2}+\left( 1-\beta _{n}\right) \left[ \left\Vert y_{n}-p\right\Vert
_{M}^{2}-\alpha _{n}\left( 1-\alpha _{n}\right) \left\Vert J_{\lambda
,M}^{A,B}\left( y_{n}\right) -y_{n}\right\Vert _{M}^{2}\right] \\
&\leq &\beta _{n}\left\Vert f\left( z_{n}\right) -p\right\Vert
_{M}^{2}+\left( 1-\beta _{n}\right) \left[ \left\Vert x_{n}-p\right\Vert
_{M}^{2}+2\theta _{n}\left\Vert x_{n}-p\right\Vert _{M}\left\Vert \text{ }%
x_{n}-x_{n-1}\right\Vert _{M}\right. \\
&&\left. +\theta _{n}^{2}\left\Vert \text{ }x_{n}-x_{n-1}\right\Vert
_{M}^{2}-\alpha _{n}\left( 1-\alpha _{n}\right) \left\Vert J_{\lambda
,M}^{A,B}\left( y_{n}\right) -y_{n}\right\Vert _{M}^{2}\right] \\
&=&\beta _{n}\left\Vert f\left( z_{n}\right) -p\right\Vert _{M}^{2}+\left(
1-\beta _{n}\right) \left\Vert x_{n}-p\right\Vert _{M}^{2}+2\left( 1-\beta
_{n}\right) \theta _{n}\left\Vert x_{n}-p\right\Vert _{M}\left\Vert \text{ }%
x_{n}-x_{n-1}\right\Vert _{M} \\
&&\left( 1-\beta _{n}\right) \theta _{n}^{2}\left\Vert \text{ }%
x_{n}-x_{n-1}\right\Vert _{M}^{2}-\alpha _{n}\left( 1-\alpha _{n}\right)
\left( 1-\beta _{n}\right) \left\Vert J_{\lambda ,M}^{A,B}\left(
y_{n}\right) -y_{n}\right\Vert _{M}^{2}
\end{eqnarray*}%
for all $n\geq 1.$ This implies that%
\begin{eqnarray*}
\alpha _{n}\left( 1-\alpha _{n}\right) \left( 1-\beta _{n}\right) \left\Vert
J_{\lambda ,M}^{A,B}\left( y_{n}\right) -y_{n}\right\Vert _{M}^{2} &\leq
&\beta _{n}\left( \left\Vert f\left( z_{n}\right) -p\right\Vert
_{M}^{2}-\left\Vert x_{n}-p\right\Vert _{M}^{2}\right) -\left\Vert
x_{n+1}-p\right\Vert _{M}^{2}+\left\Vert x_{n}-p\right\Vert _{M}^{2} \\
&&+\left( 1-\beta _{n}\right) \theta _{n}\left\Vert \text{ }%
x_{n}-x_{n-1}\right\Vert _{M}\left( 2\left\Vert x_{n}-p\right\Vert
_{M}^{2}+\theta _{n}\left\Vert \text{ }x_{n}-x_{n-1}\right\Vert _{M}\right) .
\end{eqnarray*}%
Due to the conditions $\left( iii\right) $, $\left( iv\right) $ and the
convergence of the sequence $\left\{ \left\Vert x_{n}-p\right\Vert
_{M}\right\} $, we conclude that
\begin{equation}
\lim_{n\rightarrow \infty }\left\Vert J_{\lambda ,M}^{A,B}\left(
y_{n}\right) -y_{n}\right\Vert _{M}=0.  \label{8}
\end{equation}%
On the other hand, the followings are obtained:
\begin{equation}
\lim_{n\rightarrow \infty }\left\Vert y_{n}-x_{n}\right\Vert
_{M}=\lim_{n\rightarrow \infty }\theta _{n}\left\Vert \text{ }%
x_{n}-x_{n-1}\right\Vert _{M}=0  \label{9}
\end{equation}%
and
\begin{eqnarray*}
\left\Vert z_{n}-y_{n}\right\Vert _{M} &=&\left\Vert z_{n}-J_{\lambda
,M}^{A,B}\left( y_{n}\right) +J_{\lambda ,M}^{A,B}\left( y_{n}\right)
-y_{n}\right\Vert _{M} \\
&\leq &\left\Vert z_{n}-J_{\lambda ,M}^{A,B}\left( y_{n}\right) \right\Vert
_{M}+\left\Vert J_{\lambda ,M}^{A,B}\left( y_{n}\right) -y_{n}\right\Vert
_{M} \\
&\leq &\left\Vert \left( 1-\alpha _{n}\right) y_{n}+\alpha _{n}J_{\lambda
,M}^{A,B}\left( y_{n}\right) -y_{n}\right\Vert _{M}+\left\Vert J_{\lambda
,M}^{A,B}\left( y_{n}\right) -y_{n}\right\Vert _{M} \\
&=&\left\Vert \alpha _{n}\left( J_{\lambda ,M}^{A,B}\left( y_{n}\right)
-y_{n}\right) \right\Vert _{M}+\left\Vert J_{\lambda ,M}^{A,B}\left(
y_{n}\right) -y_{n}\right\Vert _{M} \\
&=&\left( 1+\alpha _{n}\right) \left\Vert J_{\lambda ,M}^{A,B}\left(
y_{n}\right) -y_{n}\right\Vert _{M},
\end{eqnarray*}%
which implies%
\begin{equation}
\lim_{n\rightarrow \infty }\left\Vert z_{n}-y_{n}\right\Vert
_{M}=\lim_{n\rightarrow \infty }\left\Vert J_{\lambda ,M}^{A,B}\left(
y_{n}\right) -y_{n}\right\Vert _{M}=0.  \label{10}
\end{equation}%
By using (\ref{8}), (\ref{9}), (\ref{10}) and the condition $\left(
iv\right) $ we can see%
\begin{eqnarray*}
\left\Vert x_{n+1}-y_{n}\right\Vert _{M} &=&\left\Vert x_{n+1}-J_{\lambda
,M}^{A,B}\left( y_{n}\right) +J_{\lambda ,M}^{A,B}\left( y_{n}\right)
-y_{n}\right\Vert _{M} \\
&\leq &\left\Vert x_{n+1}-J_{\lambda ,M}^{A,B}\left( y_{n}\right)
\right\Vert _{M}+\left\Vert J_{\lambda ,M}^{A,B}\left( y_{n}\right)
-y_{n}\right\Vert _{M} \\
&=&\left\Vert \beta _{n}f\left( z_{n}\right) +\left( 1-\beta _{n}\right)
J_{\lambda ,M}^{A,B}\left( z_{n}\right) -J_{\lambda ,M}^{A,B}\left(
y_{n}\right) \right\Vert _{M}+\left\Vert J_{\lambda ,M}^{A,B}\left(
y_{n}\right) -y_{n}\right\Vert _{M} \\
&\leq &\beta _{n}\left\Vert f\left( z_{n}\right) -J_{\lambda ,M}^{A,B}\left(
z_{n}\right) \right\Vert _{M}+\left\Vert J_{\lambda ,M}^{A,B}\left(
z_{n}\right) -J_{\lambda ,M}^{A,B}\left( y_{n}\right) \right\Vert
_{M}+\left\Vert J_{\lambda ,M}^{A,B}\left( y_{n}\right) -y_{n}\right\Vert
_{M} \\
&\leq &\beta _{n}\left\Vert f\left( z_{n}\right) -J_{\lambda ,M}^{A,B}\left(
z_{n}\right) \right\Vert _{M}+\left\Vert z_{n}-y_{n}\right\Vert
_{M}+\left\Vert J_{\lambda ,M}^{A,B}\left( y_{n}\right) -y_{n}\right\Vert
_{M}
\end{eqnarray*}%
which implies
\begin{equation}
\lim_{n\rightarrow \infty }\left\Vert x_{n+1}-y_{n}\right\Vert _{M}=0.
\label{11}
\end{equation}

So, from the inequalities (\ref{9}) and (\ref{11}), we have%
\begin{equation*}
\left\Vert x_{n+1}-x_{n}\right\Vert _{M}\leq \left\Vert
x_{n+1}-y_{n}\right\Vert _{M}+\left\Vert y_{n}-x_{n}\right\Vert _{M}
\end{equation*}%
\begin{equation*}
\lim_{n\rightarrow \infty }\left\Vert x_{n+1}-x_{n}\right\Vert _{M}=0.
\end{equation*}

Now, we get
\begin{equation*}
\limsup_{n\rightarrow \infty }\left\langle f\left( p\right)
-p,x_{n+1}-p\right\rangle _{M}=t.
\end{equation*}

Since the sequence $\left\{ x_{n}\right\} $ is bounded, there exists a
subsequence $\left\{ x_{n_{i}}\right\} $ of $\left\{ x_{n}\right\} $ such
that $x_{n_{i}}\rightharpoonup v$ and $\lim_{i\rightarrow \infty
}\left\langle f\left( p\right) -p,x_{n_{i}+1}-p\right\rangle _{M}=t.$

By using (\ref{8}) and (\ref{9}) we can write%
\begin{eqnarray*}
\left\Vert J_{\lambda ,M}^{A,B}\left( x_{n}\right) -x_{n}\right\Vert _{M}
&=&\left\Vert J_{\lambda ,M}^{A,B}\left( x_{n}\right)
-x_{n}+y_{n}-y_{n}+J_{\lambda ,M}^{A,B}\left( y_{n}\right) -J_{\lambda
,M}^{A,B}\left( y_{n}\right) \right\Vert _{M} \\
&\leq &2\left\Vert y_{n}-x_{n}\right\Vert _{M}+\left\Vert J_{\lambda
,M}^{A,B}\left( y_{n}\right) -y_{n}\right\Vert _{M}.
\end{eqnarray*}
This implies that%
\begin{equation*}
\lim_{n\rightarrow \infty }\left\Vert J_{\lambda ,M}^{A,B}\left(
x_{n}\right) -x_{n}\right\Vert _{M}=0.
\end{equation*}

In this case, it is clear from Lemma \ref{dort} that $v\in F\left(
J_{\lambda ,M}^{A,B}\right) .$ On the other hand, since $\left\Vert
x_{n+1}-x_{n}\right\Vert _{M}\rightarrow 0$ as $n\rightarrow \infty $ and $%
x_{n_{i}}\rightharpoonup v$, we have $x_{n_{i+1}}\rightarrow v.$ Moreover,
by combining $p=P_{\Omega }f\left( p\right) $ and property of the metric
projection\ operators we can get
\begin{equation*}
\lim_{i\rightarrow \infty }\left\langle f\left( p\right)
-p,x_{n_{i}+1}-p\right\rangle _{M}=\left\langle f\left( p\right)
-p,v-p\right\rangle _{M}\leq 0.
\end{equation*}

Then this implies that
\begin{equation}
\limsup_{n\rightarrow \infty }\left\langle f\left( p\right)
-p,x_{n+1}-p\right\rangle _{M}\leq 0.  \label{12}
\end{equation}

It follows from (\ref{12}) that $\limsup_{n\rightarrow \infty }t_{n}\leq 0.$
As a result, we obtain that $x_{n}\rightarrow p.$

Secondly, we assume that there exists $n_{0}\in
\mathbb{N}
$ such that $\left\{ \left\Vert x_{n}-p\right\Vert _{M}\right\} _{n\geq
n_{0}}$ is a monotone decreasing sequence. Let us denote $\Phi_{n}
=\left\Vert x_{n}-p\right\Vert _{M}^{2}$ \ for all $n\geq 1.$ For this
reason, there exists a subsequence $\left\{ \Phi _{j}\right\} $ of $\left\{
\Phi _{n}\right\} $ such that $\Phi _{n_{j}}<\Phi _{n_{j+1}}$ for all $n\geq
n_{0}.$ Define $\tau :\left\{ n:n\geq n_{0}\right\} \rightarrow
\mathbb{N}
$ by
\begin{equation*}
\tau \left( n\right) =\max \left\{ l\in
\mathbb{N}
:l\leq n,\Phi _{l}\leq \Phi _{l+1}\right\} .
\end{equation*}

It is clear that the sequence $\tau $ is nondecreasing. By Lemma \ref{uc} we
say that $\Phi _{\tau \left( n\right) }\leq \Phi _{\tau \left( n\right) +1%
\text{ }}$ for all $n\geq n_{0.}$ So, we have
\begin{equation*}
\left\Vert \Phi _{\tau \left( n\right) \text{ }}-p\right\Vert _{M}\leq
\left\Vert \Phi _{\tau \left( n\right) +1\text{ }}-p\right\Vert _{M}.
\end{equation*}%
In the first case, by taking $\tau \left( n\right) $ instead of $n$, we can
obtain similar results. Namely, we get%
\begin{equation*}
\limsup_{n\rightarrow \infty }\left\Vert \Phi _{\tau \left( n\right) \text{ }%
}-p\right\Vert _{M}^{2}\leq 0.
\end{equation*}

Also, we have%
\begin{equation}
\left\Vert \Phi _{\tau \left( n\right) \text{ }}-p\right\Vert
_{M}^{2}\rightarrow 0\text{ and}\left\Vert \Phi _{\tau \left( n\right) +1%
\text{ }}-p\right\Vert _{M}\rightarrow 0\text{ as }n\rightarrow \infty .
\label{13}
\end{equation}

So, by using (\ref{13}) and Lemma \ref{uc}, we conclude that
\begin{equation*}
\left\Vert \Phi _{n\text{ }}-p\right\Vert _{M}\leq \left\Vert \Phi _{\tau
\left( n\right) +1\text{ }}-p\right\Vert _{M}\rightarrow 0\text{ as }%
n\rightarrow \infty .
\end{equation*}

Hence, we obtain that $x_{n}\rightarrow p,$ and the proof is completed.
\end{proof}
\end{theorem}

\section{Application to Convex Minimization Problem}

\ Now, we consider the following convex minimization problem given as a sum
of two convex functions:%
\begin{equation}
h\left( x^{\ast }\right) +g\left( x^{\ast }\right) =\min_{x\in H}\left\{
h\left( x\right) +g\left( x\right) \right\}  \label{14}
\end{equation}
Let $h:H\rightarrow
\mathbb{R}
$ is differentiable with $L_{h}$-Lipschitz gradient which is Lipschitz
constant of $\nabla h$. If $\nabla h$ is $L_{h}$-Lipschitz continuous, then
Baillon-Hadded Theorem states that\ $\nabla h$ is cocoercive with respect
to\ $L_{h}^{-1}$. Furthermore, if $g:H\rightarrow
\mathbb{R}
$ is a proper convex and lower semi-continuous function then $\partial g$ is
maximal monotone see, for detail \cite{bau}. A point $x^{\ast }$ is a
solution of minimization problem (\ref{14}) if and only if $0\in \nabla
h\left( x^{\ast }\right) +\partial g\left( x^{\ast }\right) .$ Then for any $%
\lambda >0$ we have%
\begin{eqnarray*}
0 &\in &\lambda \nabla h\left( x^{\ast }\right) +\lambda \partial g\left(
x^{\ast }\right) \\
&\Leftrightarrow &0\in \lambda L_{h}^{-1}\nabla h\left( x^{\ast }\right)
+\lambda L_{h}^{-1}\partial g\left( x^{\ast }\right) \\
&\Leftrightarrow &x^{\ast }-\lambda L_{h}^{-1}\nabla h\left( x^{\ast
}\right) \in x^{\ast }+\lambda L_{h}^{-1}\partial g\left( x^{\ast }\right) \\
&\Leftrightarrow &x^{\ast }=\left( I+\lambda L_{h}^{-1}\partial g\right)
^{-1}\left( I-\lambda L_{h}^{-1}\nabla h\right) \left( x^{\ast }\right) .
\end{eqnarray*}%
In Theorem \ref{theorem}, set $A=\partial g$, $B=\nabla h$ and $M\left(
x\right) =L_{h}x$. As a result, we can deduce the following corollary.

\begin{corollary}
Let $h:H\rightarrow
\mathbb{R}
$ be a differentiable and convex function with $L_{h}$-Lipschitz gradient
and $g:H\rightarrow
\mathbb{R}
$ be a proper convex and lower semi-continuous function. Assume that the
solution set of convex minimization problem (\ref{14}) is nonempty. The
parameter $\left\{ \theta _{n}\right\} \subset \left[ 0,\theta \right] $ and
$\left\{ \alpha _{n}\right\} ,\left\{ \beta _{n}\right\} \in \left(
0,1\right) $ satisfy the same condition as in Theorem \ref{theorem}. Let $%
\left\{ x_{n}\right\} $ be a sequence generated by%
\begin{equation}
\left\{
\begin{array}{l}
x_{0},x_{1}\in H \\
y_{n}=x_{n}+\theta _{n}\text{ }\left( x_{n}-x_{n-1}\right) \\
z_{n}=\left( I+\lambda L_{h}^{-1}\partial g\right) ^{-1}\left( I-\lambda
L_{h}^{-1}\nabla h\right) \left( \left( 1-\alpha _{n}\right) y_{n}+\alpha
_{n}\left( I+\lambda L_{h}^{-1}\partial g\right) ^{-1}\left( I-\lambda
L_{h}^{-1}\nabla h\right) y_{n}\right) \\
x_{n+1}=\beta _{n}f\left( z_{n}\right) +\left( 1-\beta _{n}\right) \left(
I+\lambda L_{h}^{-1}\partial g\right) ^{-1}\left( I-\lambda L_{h}^{-1}\nabla
h\right) z_{n}.%
\end{array}%
\right.  \label{20}
\end{equation}%
Then $\left\{ x_{n}\right\} $ converges strongly to a $x^{\ast }$solution of
convex minimization problem.
\end{corollary}

\section{Applications to Image Restoration Problem}

This section aims to show the application of the new preconditioning
forward-backward algorithm to the image restoration problem. In addition, we
conduct a comparison of the Algorithm (\ref{20}) with Algorithm (\ref{104})
and Algorithm (\ref{103}).

The inverse problem of the following form can be used to define a general
image restoration problem:
\begin{equation}
b=Ax+v  \label{16}
\end{equation}%
where $x\in
\mathbb{R}
^{d}$ is original image, $A:%
\mathbb{R}
^{d}\rightarrow
\mathbb{R}
^{m}$ is a linear operator, $b\in
\mathbb{R}
^{m}$ is observed image and $v$ is the additive noise. It is well known that
the problem (\ref{16}) is roughly comparable to a number of different
optimization problems. Also, the $l_{1}$-norm is commonly used as a
regularization tool to solve these types of problems. As a result, the image
restoration problem (\ref{16}) may be reduced to a $l_{1}$-regularization
problem, which can be expressed as
\begin{equation}
\min_{x\in
\mathbb{R}
^{d}}\left\{ \frac{1}{2}\left\Vert Ax-b\right\Vert ^{2}+\rho \left\Vert
x\right\Vert _{1}\right\} .  \label{15}
\end{equation}%
where $\rho >0$ is a regularization parameter. On the other hand, For $%
h\left( x\right) =\frac{1}{2}\left\Vert Ax-b\right\Vert ^{2}$ and $g\left(
x\right) =\rho \left\Vert x\right\Vert _{1}$, the convex minimization
problem can be reduced to $l_{1}$- regularization problem. According to this
selection, the Lipschitz gradient of $h$ is the following form $\nabla
h\left( x\right) =A^{T}\left( Ax-b\right),$ where $A^{T}$ is the transpose
of $A$.

\bigskip Now, we show that Algorithm (\ref{20}) is used to solve the image
restoration problem (\ref{16}) and also that this algorithm is compared to
Algorithm (\ref{104}) and Algorithm (\ref{103}). In all comparison, we
consider the motion and gaussian blur functions and add random noise to the
test images cameraman and mountain. In order to measure the quality of the
restored images, we use the signal to noise ratio (SNR) which is defined by%
\begin{equation*}
SNR=20\log \frac{\left\Vert x\right\Vert _{2}}{\left\Vert x-x_{n}\right\Vert
_{2}}
\end{equation*}%
where $x$ and $x_{n}$ are the original image and the estimated image at
iteration $n$, respectively. All algorithms are implemented \ in MATLAB
R2020a running on a Dell with Intel (R) Core (TM) i5 CPU and $8$ GB of RAM.

\begin{figure}[htbp]
\begin{center}
\includegraphics[width=16.0cm]{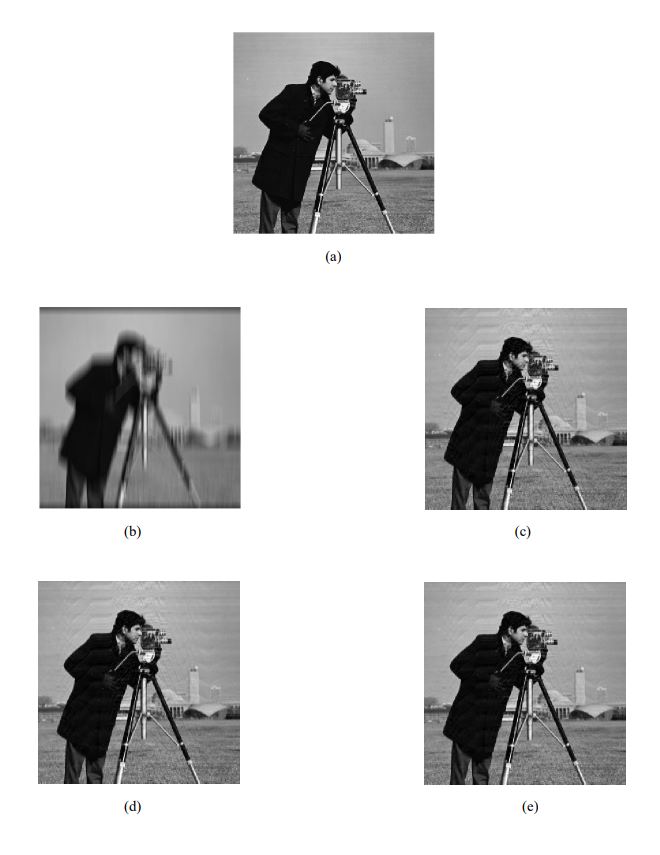}
\end{center}
\caption{ (a) Cameraman image (b) Degraded image (c) Algorithm (\protect\ref%
{103}) (d) Algorithm (\protect\ref{104}) (e) Algorithm (\protect\ref{20})}
\label{fig9b}
\end{figure}

First of all, by using cameraman image and motion blur function, we compare
Algorithm (\ref{20}) with Algorithm (\ref{104}) and Algorithm (\ref{103}).
We set $\alpha _{n}=\frac{1}{2}$, $\theta _{n}=\frac{1}{10}$, $\beta _{n}=%
\frac{1}{10n}$, $\lambda =0.99$, $f\left( x\right) =0.99x$ and the
regularization parameter $\rho=0.0001.$ Figure \ref{fig9b}, Figure \ref%
{fig10b} and Table \ref{table4} provide the visual and numerical results
corresponding to these selections.

\begin{figure}[htbp]
\begin{center}
\includegraphics[width=10.0cm]{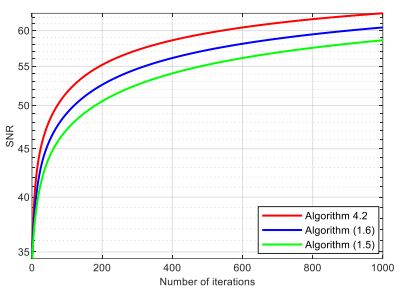}
\end{center}
\caption{ Graphic of SNR values for the Cameraman image}
\label{fig10b}
\end{figure}

\begin{table}[htbp]
\begin{center}
\begin{tabular}{cccc}
\hline
$No.Iterations$ & $Algorithm$ $(4.2)$ & $Algorithm$ $(1.6)$ & $Algorithm$ $%
(1.5)$ \\ \hline
$1$ & $35.358278$ & $34.805570$ & $34.447978$ \\
$5$ & $39.041491$ & $37.647298$ & $36.739191$ \\
$10$ & $41.596885$ & $39.737838$ & $38.459428$ \\
$25$ & $45.483306$ & $43.327179$ & $41.672360$ \\
$50$ & $48.676063$ & $46.328902$ & $44.557987$ \\
$100$ & $52.156726$ & $49.590691$ & $47.648439$ \\
$250$ & $56.836376$ & $54.268904$ & $52.177390$ \\
$500$ & $60.117373$ & $57.738798$ & $55.736154$ \\
$1000$ & $63.000553$ & $60.935109$ & $59.103851$ \\ \hline
\end{tabular}%
\end{center}
\caption{SNR values for the Cameraman image}
\label{table4}
\end{table}

Now, using mountain image and gaussian blur function, we compare Algorithm (%
\ref{20}) with Algorithm (\ref{104}) and Algorithm (\ref{103}). We take $%
\alpha _{n}=\frac{1}{2},$ $\theta _{n}=\frac{1}{2},$ $\beta _{n}=\frac{1}{2n}%
,$ $\lambda =0.99,$ and $f\left( x\right) =0.9999x.$ The numerical and
visual results corresponding to these selections are shown in Figure \ref%
{fig1b}, Figure \ref{fig2b} and Table \ref{table3}.

\begin{figure}[htbp]
\begin{center}
\includegraphics[width=16.0cm]{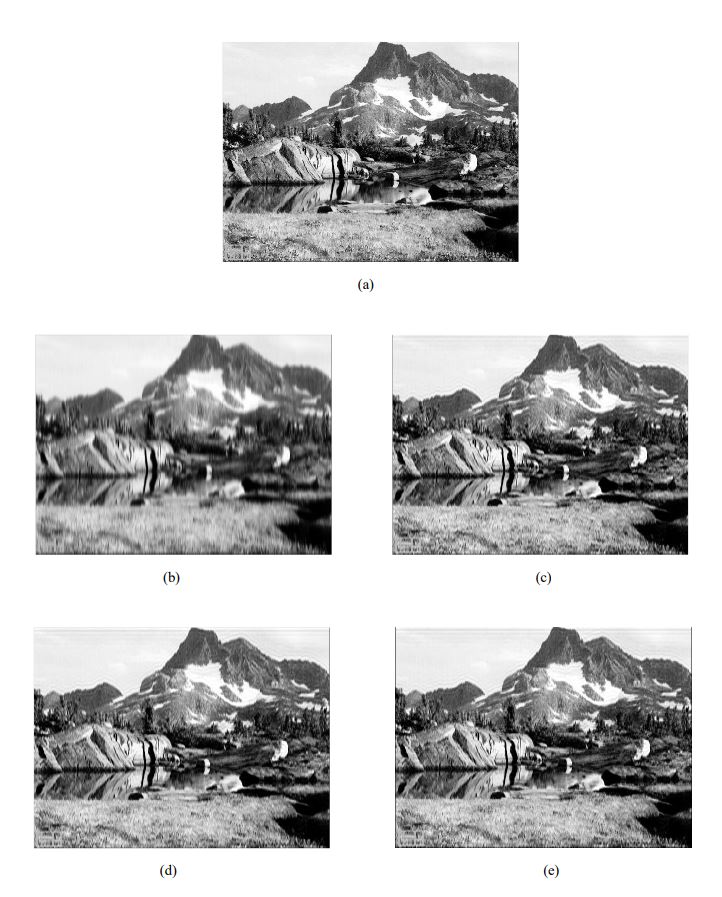}
\end{center}
\caption{ (a) Mountain image (b) Degraded image (c) Algorithm (\protect\ref%
{103}) (d) Algorithm (\protect\ref{104}) (e) Algorithm (\protect\ref{20})}
\label{fig1b}
\end{figure}

\begin{figure}[htbp]
\begin{center}
\includegraphics[width=10.0cm]{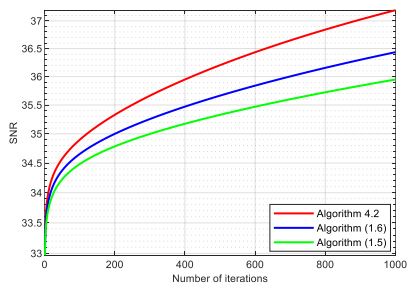}
\end{center}
\caption{ Graphic of SNR values for the Mountain image}
\label{fig2b}
\end{figure}

\begin{table}[htbp]
\begin{center}
\begin{tabular}{cccc}
\hline
$No.Iterations$ & $Algorithm$ $(4.2)$ & $Algorithm$ $(1.6)$ & $Algorithm$ $%
(1.5)$ \\ \hline
$1$ & $33.150494$ & $33.079983$ & $32.970274$ \\
$5$ & $33.975156$ & $33.862282$ & $33.758195$ \\
$10$ & $34.235762$ & $34.094985$ & $33.975192$ \\
$25$ & $34.593426$ & $34.411929$ & $34.271498$ \\
$50$ & $34.919419$ & $34.688625$ & $34.524013$ \\
$100$ & $35.350730$ & $35.031797$ & $34.822088$ \\
$250$ & $36.222684$ & $35.692589$ & $35.362060$ \\
$500$ & $37.213319$ & $36.460288$ & $35.977868$ \\
$1000$ & $38.502561$ & $37.530003$ & $36.867338$ \\ \hline
\end{tabular}%
\end{center}
\caption{SNR values for the Mountain image}
\label{table3}
\end{table}

\section{Conclusion}

In this study, we suggested a preconditioning forward-backward algorithm
which generalize some existed algorithms to handle the image restoration
problem effectively. In addition, while the weak convergence theorems were
proved for the other algorithms we generalized, we demonstrated the strong
convergence theorems for our algorithm. Experimental results demonstrate
that Algorithm (\ref{20}) restores images with a greater SNR than Algorithm (%
\ref{103}) and Algorithm (\ref{104}), indicating that its image restoration
performance is superior to Algorithm (\ref{103}) and Algorithm (\ref{104}).

\end{document}